\documentclass[12pt]{amsart}
\usepackage{latexsym,amsmath,amssymb,esint,enumitem}
\usepackage{mathtools,slashed}
\usepackage{graphicx,subfigure,psfrag}
\usepackage{hyperref}
\usepackage{xcolor}
\usepackage{cancel}
\hypersetup{
    colorlinks,
    linkcolor={red!50!black},
    citecolor={blue!50!black},
    urlcolor={blue!80!black}
}

\usepackage{graphicx}
\newcommand{\eps}{\varepsilon}
\newcommand{\RR}{{\mathbb R}}
\renewcommand{\SS}{{\mathbb S}}

\newcommand{\Vol}{\mathrm{vol}\,}

\newtheorem{theorem}{Theorem}
\newtheorem*{theorem*}{Theorem}
\newtheorem{lemma}[theorem]{Lemma}

\theoremstyle{definition}

\numberwithin{theorem}{section}
\numberwithin{equation}{section}

\begin{document}

%\setcounter{page} 0
%\tableofcontents
%\newpage

\title[A new characterization of the round sphere]
{A new characterization\\ of the round sphere}

\author[A. Burchard]{Almut Burchard} \address{University of Toronto, Canad\'a} \email{almut@math.toronto.edu}

\author[A. D. Mart\'inez]{\'Angel D. Mart\'inez} \address{CUNEF Universidad, Madrid, Spain} \email{angeld.martinezm@cunef.edu}
\date{November 5, 2026}

\begin{abstract}
We provide an elementary proof of Hersch's inequality and 
Cheng's theorem, both characterizing the sphere. We revisit these 
classical results in connection with the geodesic monotonicity 
of heat kernels on two-dimensional spheres,
which provides a third characterization.
\end{abstract}

\maketitle

%%%%%%%%%%%%%%%%%%%%%%%%%%%%%%%%%%%%%%%%%%%%%%%%%%%%%%%%%%%%%%%%%%%%%%
\section{Introduction}\label{sec-intro}
%%%%%%%%%%%%%%%%%%%%%%%%%%%%%%%%%%%%%%%%%%%%%%%%%%%%%%%%%%%%%%%%%%%%%

This note concerns three properties, two classical and one new,
each of which uniquely identifies the round metric 
on the two-dimensional sphere.  The first
characterizes the sphere as the extremizer of a 
remarkable spectral inequality,
due to Hersch~\cite{He1,He2}. It is usually stated as follows:
for any smooth metric $g$ on $\SS^2$
\begin{align}
\label{eq-Hersch}
\lambda_1(g)\operatorname{vol}_g(\mathbb{S}^2)\leq8\pi\,,
\end{align}
where $\lambda_1(g)$ is the principal eigenvalue of the 
Laplace-Beltrami operator associated with the metric $g$. 
Equality holds if and only if the metric is the standard one. It is not difficult to check that the standard metric satisfies the identity. The {\em only if} part of the later statement is more delicate and, to the best of our knowledge, it was not included 
in his original work (cf. \cite{Ch,liyau,YY}).

In any case, one arrives at the following statement.

\begin{theorem}
\label{thm-Hersch}
If equality holds in \eqref{eq-Hersch} above,
\[\lambda_1(g)\,\operatorname{vol}_g(\mathbb{S}^2) \ =\ 8\pi,\]
then $g$ is a constant multiple of the standard metric.
\end{theorem}

Digging in the literature we have been able to
cobble  together the following dispersed argument:  
Bourguignon \cite{Bourg} sketches the proof of the 
inequality and refers to Berger  \cite[p. 145 Remark 4.28]{berger}, 
for the only if part. Upon consulting the latter one learns that it 
follows from a statement in Kobayashi and Nomizu's treatise 
\cite[Example 3 on p. 341]{KN} that 
goes back to the work of Takahashi~\cite{Takahashi}. An alternative proof can be based
on results of Nadirashvili \cite{NN} (cf. \cite{elsoufi}), 
who showed that the maximizers in equation ~\eqref{eq-Hersch}
arise from minimal immersions. 
These minimal immersions are characterized
in the celebrated work of Montiel and 
Ros \cite{montielros}, which appeared more 
than a decade after Hersch's paper. 

Notice that both proofs rely on the theory of minimal surfaces and are 
not elementary in that sense. Perhaps that is the reason why
Hersch's Theorem is mentioned, but not proved,
in books intended for a general audience (such as Chavel's classic treatise \cite{Ch}). In this note we provide a self-contained proof of 
Theorem~\ref{thm-Hersch}.
%, including a proof of Hersch's inequality~\eqref{eq-Hersch}.
The proof amounts to a direct calculation with the coarea 
formula~\cite{coarea}, which is justified by 
well-known results on general elliptic equations 
with H\"older continuous coefficients \cite{Bers} and 
on the geometry of nodal sets of first eigenfunctions ~\cite{Cheng}.
It is accessible to a reader with general knowledge on Partial Differential Equations, 
even without a strong footing in the world 
of minimal immersions.  We have made an effort so that 
the presentation can be read by a strong undergraduate. We believe in this sense it also makes a valuable contribution. 

\smallskip The second property characterizes the round
sphere in terms of the first eigenspace.

\begin{theorem}[Cheng, 1976 \cite{Chengpams}]
\label{thm-Cheng}
Let $g$ be a smooth Riemannian metric on $\SS^2$.
If $\Delta_g$ has three
first eigenfunctions whose square sum is constant, 
then $g$ is a constant multiple of the 
standard metric.
\end{theorem}

The property is known to hold on the standard sphere because a basis for the first eigenspace are the spherical harmonics $x_1$, $x_2$ and $x_3$, which also are the coordinates that certainly satisfy $x_1^2+x_2^2+x_3^1=1$. We emphasize that Cheng's proof crucially relies on 
Theorem~\ref{thm-Hersch}.  
We will reverse the logic
and obtain Theorem~\ref{thm-Hersch} by
proving Theorem~\ref{thm-Cheng} first.
Both Cheng's and our argument make use of Hersch's inequality~\eqref{eq-Hersch}.

\smallskip
The third property is geodesic monotonicity 
of the heat kernel generated by the Laplace-Beltrami operator $-\Delta_g$.
In fact, our initial motivation for
this work was to classify
manifolds where the heat kernel $K_t(x,y)$
decreases as $y$ moves away from the source $x$ along a minimal geodesic. 
On the standard sphere (as on Euclidean and hyperbolic spaces) the heat kernel has the stronger property that for fixed $t$ is a decreasing function of geodesic distance alone.
We refer the reader to \cite{ACMM} 
for a description of the problem and a 
comprehensive bibliography. The basic insight is 
that manifolds with this property
are extremely rare, since monotonicity places 
strong restrictions on the eigenvalues and
eigenvectors of $-\Delta_g$~\cite{BM}.
For instance, on a flat two-dimensional torus generated by a planar lattice, 
the heat kernel is geodesically decreasing
if and only if the lattice is either rectangular
or honeycomb~\cite[Theorem 2.1]{BM}. 
Here we study metrics with monotonicity
on surfaces diffeomorphic to $\SS^2$.

\smallskip
Consider the special case of an ellipsoid $\mathcal{E}_{a_1,a_2,a_3}\subset\RR^3$ with semi-axes 
$a_1,a_2, a_3$ close to one, given by
$$\mathcal{E}_{a_1,a_2,a_3}=\left\{x\in \RR^3\ \Bigg\vert\ 
\left(\frac{x_1}{a_1} \right)^2 + \left(\frac{x_2}{a_2} \right)^2
+ \left(\frac{x_3}{a_3} \right)^2=1\right\}\,.
$$
Eigenfunctions of the Laplacian on tri-axial ellipsoids were  studied by Lam\'e~\cite{L}, among others, in the context of the heat equation. 
By a recent result of Eswarathasan and 
Kolokolnikov~\cite[equations  (30)-40)]{EK}
on perturbations of the unit sphere,
if $a_i=1+\eps\alpha_i$ for $i=1,2,3$ where
$\alpha_1\ge\alpha_2\ge\alpha_3$ and $\alpha_1+\alpha_2+\alpha_3=0$, 
then as $\eps\to 0$ the eigenvalues satisfy
\begin{align*}
\lambda_i = 2-\tfrac85 \eps\alpha_i+ O(\eps^2)
\end{align*}
for any $i=1,2,3$. For small enough $\epsilon$ least one of the first three eigenvalues is simple
unless the perturbation is trivial. It follows from a previous 
result~\cite[Theorem 2.4]{BM} that the heat 
kernel on $\mathcal{E}_{a_1,a_2,a_3}$ is not monotone on ellipsoids
near the standard sphere unless $\epsilon=0$. In particular, the heat kernel on any ellipsoid with at least 
two different semi-axes is not monotone. Suprisingly, Cheng's theorem %together with~\cite[Theorem 2.4]{BM}
allows to improve this observation beyond the perturbative regime 
to {\em any} Riemannian metric:

\begin{theorem}
\label{thm-mono} Let $g$ be a smooth Riemannian metric on the 
sphere. If there exists a sequence
of times $t_k\to\infty$ such that
the heat kernel defined by $e^{t_k \Delta}$
decreases along geodesics, then $g$ is a constant 
multiple of the standard  metric.
\end{theorem}

This provides a third characterization of the standard sphere.

\smallskip
Let us finally comment a bit on the regularity hypothesis. We have so far assumed the metric is {\em smooth}. This requirement arises because one of the first 
steps in Hersch's proof of inequality \ref{eq-Hersch} is to assume that $g$ is conformal to the 
standard metric. That assumption comes without
loss of generality if the metric is $C^{\alpha}$ a priori~\cite{Chern}.
Nevertheless, since the eigenvalues are Lipschitz continuous 
with respect to the metric (cf. \cite[Theorem 66]{yaiza}),
Hersch's inequality extends to the larger class
of continuous metrics by a straightforward continuity argument. In this paper, we largely ignore regularity issues
and assume that the metric is smooth enough.

Among the tools
employed in earlier proofs of Theorem~\ref{thm-Hersch},
the work of Takahashi~\cite{Takahashi} 
uses Riemannian connections, which involve 
first derivatives of the metric. The results of
Montiel and Ros~\cite{montielros} 
on minimal immersions
require the metric and the eigenfunctions to be at least $C^2$ a priori. 
These assumptions can be reduced by the use of a fixed conformal 
class, and by regularity theory if the metric is known to be 
smooth enough a priori. Cheng's proof of Theorem~\ref{thm-Cheng}
requires
$g\in C^2$ as it employs an identity involving 
curvature tensors. It also uses Theorem~\ref{thm-Hersch}.
Our proof of Theorem~\ref{thm-Cheng}
relies in part on the unique continuation of elliptic 
equations (cf. \cite{Cheng}).
The implication from Theorem~\ref{thm-Cheng} to~Theorem~\ref{thm-Hersch}
is valid only for metrics that are conformal to 
the standard one.

Regularity might be disregarded completely in view of Karpukhin 
and Stern's work~\cite[Theorems 1.3 and 1.4]{KS}. 
Indeed, they characterize extremizers
of Hersch's inequality by a min-max problem related 
to harmonic maps whose solutions are known 
to be smooth (cf. Chang, Wang and Yang \cite{CY} for a proof).
This is in contrast with the analogous problem for higher eigenvalues 
$\lambda_k$, $k\geq 2$, where the extremal metric has conical 
singularities. It is a pinched union of $k$ spheres, as recently 
proved by Karpukhin, Nadirashvili, Penskoi and Polterovich 
\cite{Nadirashvilietal}. This reference is also a good source
for more information on this type of problems and the 
extensive literature (cf. \cite{liyau} for manifolds of higher genus). It would be interesting to find an independent argument avoiding 
harmonic maps altogether to achieve the regularity 
of the extremal metrics a priori.

\subsection*{Organization of the paper.} 
For the sake of completeness, we prove
Hersch's of inequality~\eqref{eq-Hersch} in Section~\ref{sec-Hersch}, 
by paraphrasing his original arguments (cf.~\cite{He1,He2}). 
In Section \ref{sec-auxiliary} 
we provide auxiliary results from the literature
on first eigenfunctions, and derive some consequences that 
will be useful in the proof of Cheng's theorem.
In Section \ref{sec-Cheng} we prove Cheng's result by 
elementary means, as promised. 
The final Sections~\ref{sec:6} and \ref{sec:7}
contain the proofs of Theorems~\ref{thm-Hersch}
and~\ref{thm-mono}, respectively. To ease the notation, 
all geometric quantities (inner products, moduli, gradients,  
volume, arc length, and eigenvalues) are understood 
relative to the metric $g$.  
This provides some clarity (at a price). 
The dependence on the metric
will be made explicit whenever we compare $g$ to 
the standard metric on $\SS^2$.

%%%%%%%%%%%%%%%%%%%%%%%%%%%%%%%%%%%%%%%%%%%%%%%%%%%%%%%%%%%%%%%%%%%%%%
\section{Proof of Hersch's inequality}\label{sec-Hersch}
%%%%%%%%%%%%%%%%%%%%%%%%%%%%%%%%%%%%%%%%%%%%%%%%%%%%%%%%%%%%%%%%%%%%%%

Hersch's %original 
proof of equation ~\eqref{eq-Hersch}
is spread over two of his articles
\cite{He1,He2}, published with a separation of almost ten 
years. At the moment of this writing they are not accessible, 
although we have been fortunate to locate translations 
of them due to C. Shonkwiler. Here we revisit the argument. 
Another sketch of this proof can be found in \cite{N}.

The proof starts with a variational 
inequality for the inverses of the eigenvalues
of $-\Delta_g$ on a manifold $M$.
This requires a bit of linear algebra.
Let $L$ be a real linear space of dimension $n$.
Given two symmetric bilinear forms
$A$ and $B$, with $B$ positive definite,
define the {\em trace of $A$ with respect to $B$} by
\[ \operatorname{Tr}_B(A)=\sum_{i=1}^nA(\nu_i,\nu_i),\]
where $\{\nu_1,\ldots,\nu_n\}$ is a $B$-orthonormal basis
for $L$.  It is elementary to show that the trace 
does not depend on the choice of the $B$ orthonormal basis: suppose we have another base given by
\[\nu_i=\sum_{k=1}^n\beta_{ik}\mu_k\]
where $\beta_{ij}$ is an orthonormal matrix. Indeed, 
\[\delta_{ij}=B(\nu_i,\nu_j)=\sum_{k,k'=1}^n\beta_{ik}\beta_{jk'}B(\mu_k,\mu_{k'})=\sum_{k,k'=1}^n\beta_{ik}\beta_{jk'}\delta_{kk'}=\sum_{k=1}^n\beta_{ik}\beta_{jk}.\]
A straightforward calculations shows using $\beta_{ij}$ orthogonality concludes that the trace is well-defined. This will be useful in the following

\begin{lemma}[Hersch's variational
inequality~\cite{He1}]\label{lem-var}
Given a collection of~$n$ linearly independent
differentiable functions $u_1$,\ldots, $u_n$ 
on $(M,g)$ that are
orthogonal to the constants and satisfy 
\begin{equation}
\label{eq-B-ortho}
\int_M \langle \nabla u_i,\nabla u_j\rangle\,
d\operatorname{vol} =0
\end{equation}
where $1\le i<j\le n$. Then
\begin{align}
\label{eq-var}
\sum_{i=1}^n\lambda_i^{-1}\geq\sum_{i=1}^n\frac{\int_M| u_i|^2d\operatorname{vol}}
{\int_M |\nabla u_i|^2d\operatorname{vol}}
\end{align}
holds.\end{lemma}

\begin{proof}
Let $L$ be the span of $u_1,\ldots, u_n$.
By definition, the right hand side is the trace of the bilinear form
\[A(u,u)=\int_M |u|^2\, d\operatorname{vol}\]
with respect to
\[ B(u,u)=\int_M |\nabla u|^2\,d\operatorname{vol}\,.\]
Note that $A$ and $B$ are positive definite
on $L$, and the denominators in equation ~\eqref{eq-var}
serve to $B$-normalize the functions $u_1,\ldots, u_n$.

Since $\operatorname{Tr}_B(A)$ is independent of the 
choice of basis, it is enough to show 
that the inequality holds for a specific 
set of $B$-orthogonal functions. 
We will construct a $B$-orthogonal set
$v_1,\ldots, v_n$ in $L$ such that the Rayleigh quotients satisfy
\begin{align}
\label{eq-Rayleigh}
\lambda_i\ \ge\ 
\frac{\int_M |\nabla v_i|^2d\operatorname{vol}} {\int_M| v_i|^2d\operatorname{vol}}
\qquad (i=1,\ldots, n)\,.
\end{align}

To that end let $\varphi_1$,\ldots, $\varphi_n$ be $L^2$-orthogonal
eigenfunctions for the eigenvalues $\lambda_1,\ldots, \lambda_n$.
We observe that $L$, having dimension $n$,
contains a function, $v_n$, 
that is $L^2$ orthogonal to $\varphi_1,\ldots,\varphi_{n-1}$.
By the min-max characterization
of eigenvalues, equation ~\eqref{eq-Rayleigh} holds for $i=n$. Recursively, once $v_{k+1},\ldots v_n$ have been constructed, 
choose $v_k\in L$ to be 
$B$-orthogonal to $v_{k+1},\ldots v_n$
and to $\varphi_1,\ldots, \varphi_{k-1}$.
Note that $v_k$ is also $L^2$-orthogonal to 
$\varphi_1,\ldots,\varphi_{k-1}$, since these are
eigenfunctions associated with non-zero eigenvalues. 
Therefore equation ~\eqref{eq-Rayleigh} holds also for $i=k$.
Equation~\eqref{eq-var} follows by adding the reciprocals
of the inequalities in equation ~\eqref{eq-Rayleigh}.
Clearly, equality is achieved if $L$ has a basis consisting of $n$ lowest eigenfunctions.
\end{proof}

In addition to this variational inequality, the 
proof of equation ~\eqref{eq-Hersch} takes advantage of the 
Uniformization Theorem and the conformal invariance of the gradient 
square integral.
Indeed, any  smooth
metric on the two dimensional sphere can be written 
(changing coordinates) as a conformal deformation
of the standard metric. In other words, we can assume without loss of generality
that $g=e^f g_0$, where $g_0$ denotes the standard
metric on the sphere. 

Let us fix coordinates by choosing a conformal diffeomorphism
$X:\SS^2\to\SS^2$. Viewing the map $X$ into the unit sphere as an embedding into $\RR^3$, we let
$X_i$ be the $i$-th component for any $i=1,2,3$.  We will 
apply Lemma~\ref{lem-var}
to the functions $u_i=X_i$. The requisite orthogonality of the $X_i$ to constants
is equivalent to the center of mass
of the vector $X$ being the origin, i.e.,
\begin{align}
\label{eq-vanish}
\int_{\SS^2} X d\operatorname{vol}_g \ =\ (0,0,0)\,.
\end{align}
We claim that if the given embedding $X$ does not
have this property, we can achieve it
by composing $X$ it with a conformal map, as follows.

\begin{lemma}
\label{lem-vanish}
Let  $\mu$ be a  probability
measure on $\mathbb{S}^2$ that contains no point masses. 
There exists
a conformal transformation $\psi:\mathbb{S}^2\rightarrow \mathbb{S}^2$ such that
\[ \int_{}\psi \circ Xd\mu=(0,0,0).\] 
\end{lemma}

\begin{proof}  Let us define the center of mass by
\[V=\int_{}Xd\mu.\]
If the vector $V=(0,0,0)$ we can take $\psi$ the identity and there is nothing to prove. Otherwise, let us denote by $\sigma_p:
\mathbb{S}^2\rightarrow T_p\mathbb{S}^2$ the 
stereographic projection from  the sphere to the tangent 
plane at the point $p$, and by $m_{p,t}$ the multiplication by $t$ 
in the tangent plane $T_p\mathbb{S}^2$.
Finally, for $p\in\SS^2$ and $t>0$
we define the conformal maps $\psi_{p,t}=
\sigma^{-1}_p\circ m_{p,t}\circ \sigma_{p}$.

For $t$ very small $\psi_{p,t}$
has the effect of concentrating all the mass 
near $p$. (Here we are using that there are no point masses.) 
On the other hand $\psi_{p,1}$ is the identity and hence 
the center of mass is at $V$, regardless of $p$.
This shows that the map 
\[\Psi_t: p\mapsto \int_{\mathbb{S}^2}\psi_{p,t}\circ X d\mu\] 
is nearly the identity map on $\mathbb{S}^2$ provided $t$ is small  enough,
while it is nearly the constant $V\neq 0$ for $t$ close to 
one. Moreover, $|\Psi_t(p)| \le 1$ for all $t\in [0,1]$ and
$p\in\mathbb{S}^2$ by convexity.
Thus $\Psi_t$ defines a homotopy that contracts
$\mathbb{S}^2$ inside the closed unit ball to the vector $V$.
Since such a homotopy swipes our the entire
unit ball, there exist $p_0\in\mathbb{S}^2$ and $t_0\in (0,1)$ such that
$\Psi_{t_0}(p_0)=0$, as claimed.
\end{proof}

Continuing with the proof of equation ~\eqref{eq-Hersch},
we use Lemma~\ref{lem-vanish}
with $\mu=\operatorname{vol}_g$ to obtain
a conformal transformation $\psi$ such that
equation ~\eqref{eq-vanish} holds, without loss of generality. By conformal invariance of the Dirichlet integral and a change of variables,
\begin{align*}
\int_{\mathbb{S}^2}\langle\nabla_g X_i,\nabla_gX_i\rangle_g\, d\operatorname{vol}_g
=\int_{\mathbb{S}^2}\langle \nabla_{g_0}X_i,\nabla_{g_0} X_j\rangle_{g_0}d\operatorname{vol}_{g_0}
&=\frac23 \operatorname{vol}_{g_0}(\mathbb{S}^2)\delta_{ij}=\frac{8\pi}{3}\delta_{ij}
\end{align*}
where we have used integration by parts and the fact that the
$X_i$ are eigenfunctions of $-\Delta_{g_0}$
with eigenvalue $2$ in the second identity.
This verifies the orthogonality conditions in
equation~\eqref{eq-B-ortho}. By Lemma \ref{lem-var},
\begin{align}
\label{eq-Hersch-proof}
3\lambda_1^{-1} \ge \sum_{i=1}^3\lambda_i^{-1}
\geq\sum_{i=1}^3\frac{\int_{\mathbb{S}^2} |X_i|^2 e^fd\operatorname{vol}_{g_0}}
  {\int_{\mathbb{S}^2 }|\nabla X_i|^2d\operatorname{vol}_{g_0}} = \frac{3}{8\pi}\operatorname{vol}_g(\mathbb{S}^2)\,,
\end{align}
where the first inequality uses that 
$\lambda_1\le \lambda_2\le \lambda_3$, and
the last identity uses that 
$\int_{\mathbb{S}^2}|\nabla X_i|^2\, d\operatorname{vol}_g=\frac{8\pi}{3}$
and $|X|^2=1$.
Solving for $\lambda_1$ yields equation ~\eqref{eq-Hersch}. This completes the proof of Hersch's inequality.

\bigskip \noindent
\textbf{Remark:} \ Equality holds in the case of the standard sphere. 
Conversely, if equality is achieved then the $X_i$ are forced, by the min-max characterization of $\lambda_1$ and equation \eqref{eq-Hersch-proof},  to 
be first eigenfunctions of $-\Delta_g$. Furthermore, since $|X|=1$ they satisfy the hypothesis of Theorem~\ref{thm-Cheng}. This observation will be useful later.

%%%%%%%%%%%%%%%%%%%%%%%%%%%%%%%%%%%%%%%%%%%%%%%%%%%%%%%%%%%%%%%%%%%%%%
\section{Properties of first eigenfunctions}\label{sec-auxiliary}
%%%%%%%%%%%%%%%%%%%%%%%%%%%%%%%%%%%%%%%%%%%%%%%%%%%%%%%%%%%%%%%%%%%%%%

The following results of Cheng \cite{Cheng} will play a 
key r\^ole in the proof of Theorem~\ref{thm-Cheng}.

\begin{lemma}[{\cite[Theorem 3.2 and Corollary 3.5]{Cheng}}]
\label{lem-phi1}
For an arbitrary smooth Riemannian metric on $\mathbb{S}^2$,
the multiplicity of the principal eigenvalue is at most three.  Moreover, for any first eigenfunction, $0$ is a regular value, 
and the nodal set is a smooth closed curve.
\end{lemma}
This result relies on an analysis of the equation that goes back 
to the work of Bers \cite{Bers} and the existence of a 
certain diffeomorphism as in \cite{Cheng}. Pertinent to our regularity comments in the introduction, notice that the nodal set might not be even be a curve if the metric is not smooth enough, due to non-uniqueness under the Lipschitz threshold. Lemma~\ref{lem-phi1} together with the Jordan Curve Theorem implies that 
the nodal set divides the sphere into at least two connected domains 
where the eigenfunction is respectively positive and negative. 
Courant's nodal domain theorem proves that there are {\em exactly} two nodal domains. The next three lemmas investigate consequences of the hypotheses of
Theorem~\ref{thm-Cheng}.

\begin{lemma} 
\label{lem-grad-squares}
Let $\varphi_1,\ldots, \varphi_m$ be eigenfunctions
associated with the same  eigenvalue, $\lambda$. If
\begin{align}
\label{eq-squares} 
\varphi_1^2+\ldots+\varphi_m^2=c^2
\end{align}
for some constant $c$, then
\begin{equation} \label{eq-grad-squares}
|\nabla\varphi_1|^2+ \ldots |\nabla\varphi_m|^2=\lambda c^2\,.
\end{equation}
%(cf. Lemma 5.1 and Lemma 5.4 \cite{BM} respectively).
\end{lemma}

\begin{proof} By the product rule
\begin{align*}
\sum_{i=1}^m |\nabla\varphi_i|^2\ 
&= \ \sum_{i=1}^m \left(\frac12 \Delta (\varphi_i^2)- \varphi_i\Delta\varphi_i\right)\\
&= \ \frac12 \Delta\left(\sum_{i=1}^m \varphi_i^2\right) +
\sum_{i=1}^m \lambda \varphi_i^2\\
&= \ \lambda c^2\,,
\end{align*}
where equation~\eqref{eq-squares} was used in the last identity.
\end{proof}

\begin{lemma}
\label{lem-independent}
For any smooth Riemannian metric on $\mathbb{S}^2$, if 
three first eigenfunctions $\varphi_1, \varphi_2, \varphi_3$ satisfy
\begin{equation}
\label{eq-constant}
\varphi_1^2+\varphi_2^2+\varphi_3^2=c^2
\end{equation}
for some $c>0$, then they are linearly independent.
In particular, the principal eigenvalue
$\lambda_1$ has multiplicity three.
\end{lemma}

\begin{proof} We follow the line of argument 
in~\cite[Proof of Theorem 2, and Section 6]{BM}.
Let $L$ be the vector space spanned by
the three eigenfunctions. Note that $L$ is orthogonal
to the constant functions. Clearly, $\dim L>1$, since otherwise all 
functions in $L$ would have a common zero, contradicting
the assumption that $c>0$. In order to exclude $\dim L=2$ we observe that the nodal sets of any two principal eigenfunctions $\varphi_1,\varphi_2$ has non-empty intersection. Otherwise a nodal domain of $\Omega$ of one of the eigenfunctions (say, $\varphi_1$) would be compactly contained in a nodal domain of the other. But then, in Green's identity 
\[
\int_{\Omega}
\left(
\varphi_1 \Delta \varphi_2
-
\varphi_2 \Delta \varphi_1
\right)\,d\mathrm{vol}
=
\int_{\partial\Omega}
\left(
\varphi_1 \frac{\partial \varphi_2}{\partial \nu}
-
\varphi_2 \frac{\partial \varphi_1}{\partial \nu}
\right)\,d\sigma
\]
the left hand side vanishes for eigenfunctions with the same eigenvalue. On the right hand side the first integral vanishes while the second has a sign (the normal derivative has a sign due to Hopf's lemma), a contradiction. It follows that necessarily $\dim L=3$, which means that the three eigenfunctions
are linearly independent and span the first eigenspace.
\end{proof}

Under the assumptions of Lemma~\ref{lem-independent},
let $\Phi:\mathbb{S}^2\to\mathbb{S}^2$ be the map with
components $c^{-1}\varphi_i$. 
We will use $\Phi$ to define coordinates
on $\mathbb{S}^2$. For that purpose it is enough to notice the following 

\begin{lemma} \label{lem-diffeo}
The map $\Phi$ is a diffeomorphism.
\end{lemma}

\begin{proof}  We need to prove that $\Phi$ is bijective,
and that its derivative $D\Phi(p)$ is invertible 
for every $p\in\mathbb{S}^2$. The proof proceeds in three steps.

\smallskip
(1) {\em $D\Phi$ has full rank.} \
Suppose that at some point $p\in\mathbb{S}^2$ the derivative $D\Phi(p)$ 
has rank less than two.  By composing $\Phi$ with a suitable rotation 
in the image, we may assume that
$\Phi(p)$ is the north pole of $\mathbb{S}^2$. Let $u,v,w$ denote 
the components of $\Phi$. Then $u(p)=v(p)=0$, and $w(p)=1$.
Since $D\Phi(p)$ has rank less than two and $\nabla w(p)=0$, 
the gradients $\nabla u(p)$ and $\nabla v(p)$ are linearly dependent.
Therefore we can find a nontrivial linear combination
$$a\nabla u(p)+b \nabla v(p)=0.$$
The eigenfunction $\varphi= a u+b v$ has $p$ in its 
nodal set, and its gradient vanishes at $p$. 
But this contradicts Lemma~\ref{lem-phi1}.
Therefore $\nabla u(p)$ and $\nabla v(p)$ are linearly independent, 
and $D\Phi(p)$ has rank two. 

\smallskip
(2) {\em $\Phi$ is surjective.} 
By (1) and the Inverse Function Theorem, 
$\Phi$ is a local diffeomorphism. In particular its range 
is relatively open in $\mathbb{S}^2$. Since its domain is compact, 
its range is also closed. It follows by  connectedness 
of $\mathbb{S}^2$ that $\Phi$ is surjective.

\smallskip 
(3) {\em $\Phi$ is injective.} \
For any point $\xi\in\mathbb{S}^2$, the
inverse image $\Phi^{-1}(\xi)$ is a compact subset of $\mathbb{S}^2$.  Since $\Phi$ 
is a local diffeomorphism, this subset is finite,
and $\Phi$ is a covering map.  Since the base space
$\mathbb{S}^2$ is simply connected, the covering has  a
single leaf.
Therefore $\Phi$ is a diffeomorphism.

\end{proof}

Finally, we show that each of the coordinate functions
provided by the components of $\Phi$ has exactly two critical points.

\begin{lemma} 
\label{lem-Morse}
Under the assumptions of Lemma~\ref{lem-independent},
the gradient of any first eigenfunction vanishes at exactly two points,
where it attains its maximum and minimum.
\end{lemma}

\begin{proof} Let $\varphi$ be a first eigenfunction, and let $p\in \mathbb{S}^2$ 
be a point where $\nabla \varphi$ vanishes.  By composing the diffeomorphism $\Phi$ with a rotation, we may take 
$\Phi(p)$ to be the north pole, $(0,0,1)$.
Denoting the components of $\Phi$ by $u$, $v$, and $w$, we express 
the eigenfunction as a linear combination,
$\varphi=au+bv+cw$. Clearly $\nabla w(p)=0$ because $w$ 
attains its maximum there. It follows that 
$$a\nabla u(p) +b\nabla v(p)=\nabla \varphi(p)=0.$$ 
Since $\nabla u(p)$ and $\nabla v(p)$ are linearly independent 
by Lemma~\ref{lem-diffeo}, we have $a=b=0$, i.e., $\varphi=cw$.
Therefore $\varphi$ attains either 
its maximum or minimum at~$p$, depending on the sign of $c$.
\end{proof}

%%%%%%%%%%%%%%%%%%%%%%%%%%%%%%%%%%%%%%%%%%%%%%%%%%%%%%%%%%%%%%%%%%%%%%
\section{Proof of Theorem \ref{thm-Cheng}} \label{sec-Cheng}
%%%%%%%%%%%%%%%%%%%%%%%%%%%%%%%%%%%%%%%%%%%%%%%%%%%%%%%%%%%%%%%%%%%%%%

We follow Cheng's basic strategy and
start by showing that under the hypothesis of Theorem~\ref{thm-Cheng},
Hersch's inequality~\eqref{eq-Hersch} is reversed.
However, our proof of this fact is new. 
Instead of Cheng's curvature tensors estimates
 (see~\cite[Theorem 1]{Chengpams}),
we apply the coarea formula
to the gradient of certain coordinate functions.

\begin{lemma} \label{lem-reverse}
Let $g$ be a smooth Riemannian metric on $\mathbb{S}^2$.
If $\Delta_g$ has three first eigenfunctions
$u,v,w$ with
$u^2+v^2+w^2=1$ on $\mathbb{S}^2$, then
\begin{align}
\label{eq-reverse}
\lambda_1(g)\,\operatorname{vol}_g(\mathbb{S}^2)\ge 8\pi\,.
\end{align}
\end{lemma}

\begin{proof} We introduce cylindrical coordinates
in the range of $\Phi$ by writing
$$
u+iv= \sqrt{1\!-w^2}\, e^{i\alpha}\,.
$$
By definition, for any point $p\in \mathbb{S}^2$
the angular variable $\alpha(p)$ is
the longitudinal coordinate of $\Phi(p)$ in $\mathbb{S}^2$,
and $w(p)$ is the latitude.  As long as  $w\ne \pm 1$, 
the value of $\alpha$ is defined up to an integer multiple of $2\pi$. 

By the chain rule
$$
\nabla u + i\nabla v = \left(-\frac{w}{\sqrt{1\!-\!w^2}}\nabla w 
+ i \sqrt{1\!-\!w^2} \nabla \alpha\right)e^{i\alpha}\,,$$ 
and for the squared modulus
$$
|\nabla u|^2 + |\nabla v|^2 = 
\frac{w^2}{1\!-\!w^2}\,|\nabla w|^2
+ (1\!-\!w^2)\,|\nabla\alpha|^2.
$$
The gradient field $\nabla\alpha$ is well-defined 
and smooth except where 
$w=\pm 1$. It follows from equation ~\eqref{eq-grad-squares} that
\begin{align}
    \label{eq-grad-alpha} 
\begin{split} \lambda_1 
& = |\nabla u|^2 +|\nabla v|^2 +|\nabla w|^2\\
& = \frac{~1\,}{1\!-\!w^2}\,|\nabla w|^2
+ (1\!-\!w^2)\,|\nabla\alpha|^2 .
\end{split}
\end{align}
By Lemma \ref{lem-Morse} the level
sets $\{w=h\}$ are smooth simple closed curves for 
all $h\in (-1,1)$.  Let
\[\ell(h):=\int_{\{w=h\}} 1\, d\sigma
\]
be the {\em length} of the curve $\{w=h\}$  .
We average equation ~\eqref{eq-grad-alpha} over the curve to obtain
\begin{align}\label{elementary*}\lambda_1 
 =
\frac{ ~1\,}{(1\!-\!h^2)\ell(h)}\,\int_{\{w=h\}}|\nabla w|^2
d\sigma + 
\frac{(1\!-\!h^2)}{\ell(h)}\,\int_{\{w=h\}}|\nabla\alpha|^2 \, 
d\sigma.\end{align}
Applying Jensen's inequality to each of the integrals, followed by the elementary inequality $\frac{a^2}{\delta}+\delta b^2\ge 2ab$, one gets
\begin{align}
\label{eq-CS} 
\lambda_1
&\ge\ 
\frac{2}{\ell(h)^2}
\left(\int_{\{w=h\}}|\nabla w| d\sigma\right)
\left(\int_{\{w=h\}}|\nabla\alpha| d\sigma\right)\,.
\end{align}
The denominator is estimated with the aid of Cauchy-Schwarz by
\begin{align*}
\ell(h)^2 \le \left(\int_{\{w=h\}} |\nabla w|\, d\sigma\right)
\left(\int_{\{w=h\}} |\nabla w|^{-1}\, d\sigma\right)\,.
\end{align*}
We insert this into 
equation ~\eqref{eq-CS}, cancel the
common factor, and collect terms to obtain
\begin{align}
\label{eq-lambda}
\lambda_1 \int_{\{w=h\}} |\nabla w|^{-1}\, d\sigma
\ge 
2 \int_{\{w=h\}} |\nabla\alpha|\, d\sigma\,.
\end{align}
The integral on the left hand side
we recognize (with the help of 
the coarea formula, see \cite{coarea}) as the derivative 
of the distribution function, 
$$
\int_{\{w=h\}}|\nabla w|^{-1}\, d\sigma
= \frac{d}{dh} \Vol(\{w\le h\})=:\mu'(h)\,.$$ 
For a proof of the coarea formula we refer to \cite[Theorem 1 in 
IV.1]{Ch}.
The differentiation is justified 
because the only critical values of $w$ are its extremes 
$\pm 1$ (cf. Lemma~\ref{lem-Morse}).

The integral on the right hand side we estimate by
\begin{align}
\label{eq-alpha}
\int_{\{w=h\}}|\nabla\alpha|\, d\mathcal{H}^1
\ \ge \ 
\left\vert\int_{\{w=h\}}
\langle \nabla\alpha,\tau\rangle \, d\mathcal{H}^1 \right\vert
\ = \ 2\pi\,,
\end{align}
where $\tau$ is the unit tangent to
the curve $\{w=h\}$. The equality on the right holds
since $\{w=h\}$
is a circle of constant latitude, along which the value of $\alpha$ 
changes by $2\pi$. 
We arrive at the inequality $\lambda_1\mu'(h)\ge 4\pi$.
Integration over $h\in [-1,1]$ yields the claim.
\end{proof}

\smallskip To conclude the proof,  we need to compare the given metric $g$ on $\SS^2$ with the standard metric $g_0$. To avoid confusion, from now on the metric for geometric quantities such as the volume will be specified by a subscript
($g$ or $g_0$).
We use the coordinates introduced in the proof of Lemma~\ref{lem-reverse} on $(\SS^2,g)$, namely the height $w=z$ and the angular variable $\alpha=\mathrm{arg}(x+iy)$. The result follows if in those coordinates, the metric
\[g=\left(\begin{array}{cc}
g_{\alpha \alpha}&g_{\alpha w}\\g_{w\alpha}&g_{ww}\end{array}\right)
\]
is a constant multiple of the standard metric
\[g_0=\left(\begin{array}{cc}
1-w^2&0\\0&(1-w^2)^{-1}\end{array}\right).\]

So far we have reversed Hersch's inequality~\eqref{eq-Hersch}
under the assumption that $|\Phi|~=~1$. Therefore
\eqref{eq-Hersch}
holds with equality. 
At this point, Cheng invokes
Theorem~\ref{thm-Hersch}. Instead, we analyze the equality 
cases of the reverse inequality~\eqref{eq-reverse}
by closely inspecting the proof of Lemma~\ref{lem-reverse}. In particular equality holds
in equations  \eqref{eq-CS}, \eqref{eq-lambda} and \eqref{eq-alpha} for every $h\in (-1,1)$.

We may assume, by scaling, that $\Vol_g(\SS^2)=\Vol_{g_0}(\SS^2)=4\pi$, and thus
$$
\lambda_1(g)=\lambda_1(g_0)=2\,.
$$
Consider first~\eqref{eq-CS}. Equality in Jensen's inequality implies that 
$|\nabla_g w|_g$ and $|\nabla_g \alpha|_g$ are equal to constants (depending on $h$) on each level curve, and equality in the subsequent elementary inequality implies
\[ 
|\nabla_g w|_g \ =  (1-h^2)|\nabla_g \alpha|\,.
\]
Since these are constant on $\{w=h\}$,
equality in \eqref{eq-lambda} further implies
$$|\nabla_g w|_g\, |\nabla_g \alpha|_g=\frac12 \lambda_1(g)=1.$$
We conclude that
$|\nabla_g \alpha|_g=\sqrt{1-h^2}=|\nabla_g w|_g^{-1}$. This determines
the diagonal entries of the cometric, which are computed with the aid of the scalar product definition as
\[
\left\{\begin{aligned}
g^{\alpha\alpha}&=\langle \nabla_g \alpha,\nabla_g\alpha\rangle_g
=(1-w^2)^{-1},\\
g^{ww}& =\langle \nabla_g w,\nabla_g w\rangle_g=1-w^2\,.
\end{aligned}\right.
\]
Moreover, equality in equation ~\eqref{eq-alpha}
implies that $\nabla_g \alpha$ is tangent
to the level curve $\{w=h\}$, and hence 
orthogonal to $\nabla_g w$. Therefore the off-diagonal
entries of the cometric vanish,
\begin{align*}
g^{\alpha w}&=\langle \nabla_g \alpha,\nabla_gw\rangle_g=0\,.
\end{align*}
In summary,
\[
g^{-1} 
\ =\ \left(\begin{array}{cc} (1-w^2)^{-1}&0\\0&1-w^2
\end{array}\right) 
\ =\ (g_0)^{-1},
\]
concluding the proof of Theorem~\ref{thm-Cheng}

%%%%%%%%%%%%%%%%%%%%%%%%%%%%%%%%%%%%%%%%%%%%%%%%%%%%%%%%%%%%%%%%%%%%%%

%%%%%%%%%%%%%%%%%%%%%%%%%%%%%%%%%%%%%%%%%%%%%%%%%%%%%%%%%%%%%%%%%%%%%%

\section{Proof of Theorem~\ref{thm-Hersch}} \label{sec:6}
Let $g$ be a smooth metric for which
equality holds in
Hersch's inequality~\eqref{eq-Hersch}. As remarked at the end of the proof the coordinates of $X$ are three first eigenfunctions. This means that the principal eigenvalue
has multiplicity three, with 
eigenfunctions
given by the coordinate functions $X_1$, $X_2$, and $X_3$
that satisfy $X_1^2+X_2^2+X_3^2=1$ by construction.
It follows from Theorem~\ref{thm-Cheng} that
the metric $g$ is a constant multiple of the standard metric, concluding the proof.

\section{Proof of Theorem \ref{thm-mono}}\label{sec:7}

We begin by recalling some auxiliary results from previous work~\cite{BM}.

\begin{lemma}[{\cite[Lemma 5.2]{BM}}]
\label{lem-squares-m}
\ Let $(M, g)$ be a compact connected Riemannian manifold 
whose principal eigenvalue has multiplicity 
$m$, and let $\varphi_1,\ldots, \varphi_m$
be an orthonormal basis of eigenfunctions for 
$\lambda_1(g)$.  If the heat kernel 
associated with $e^{t\Delta_g}$
decreases along geodesics for some sequence of 
times $t_k\to \infty$, then
\begin{equation}
\label{eq-squares-m}
\varphi_1^2 +\cdots + \varphi_m^2=\frac{m}{\operatorname{vol}_g(M)}\,.
\end{equation}
\end{lemma}

Under the assumptions of Theorem~\ref{thm-mono},
Lemma~\ref{lem-squares-m} implies that the eigenfunctions
associated with the principal eigenvalue satisfy 
equation ~\eqref{eq-squares-m}.
By Lemmas~\ref{lem-phi1} and
~\ref{lem-independent} the multiplicity
of $\lambda_1(g)$ is $m=3$, and
$\varphi_1^2+ \varphi_2^2+ \varphi_3^2$ is constant.
By Theorem~\ref{thm-Cheng} the metric $g$ is a 
constant multiple of the standard metric.

\subsection*{Acknowledgments.} 
Partial funding for the research in this
paper was provided by NSERC Discovery Grant \# 311685.
The second named author was partially supported by a 
Mar\'ia Zambrano Postdoctoral Fellowship.

\end{document}